\documentclass[a4paper,12pt,oneside,reqno]{amsart}

\pdfoutput=1

\usepackage[T1]{fontenc}
\usepackage[utf8]{inputenc}
\usepackage{lmodern}
\usepackage[english]{babel}

\usepackage[%
	left=2.5cm,       
	right=2.5cm,      
	top=3.5cm,        
	bottom=3.5cm,     
	heightrounded,    
	bindingoffset=0mm 
]{geometry}

\usepackage{amsmath}
\usepackage{amssymb}
\usepackage{mathtools}
\usepackage{mathrsfs}
\usepackage[
colorlinks=true,
urlcolor=teal,
citecolor=magenta,
linkcolor=teal,
linktocpage,
pdfpagelabels,
bookmarksnumbered,
bookmarksopen,
]
{hyperref}
\usepackage[
capitalize, 
nameinlink,
noabbrev]
{cleveref}
\usepackage{bookmark}

\usepackage[initials,
msc-links]{amsrefs}

\usepackage{xfrac}
\usepackage{graphicx}
\usepackage{tikz}
\usepackage{pgfplots}
\pgfplotsset{compat = newest}
\usepackage{xcolor}
\usepackage{enumitem}
\usepackage{url}

\numberwithin{equation}{section}

\theoremstyle{plain}
\newtheorem{theorem}{Theorem}[section]

\newtheorem{corollary}[theorem]{Corollary}

\theoremstyle{definition}

\newtheorem{remark}[theorem]{Remark}

\newcommand{\R}{\mathbb{R}}

\newcommand{\eps}{\varepsilon}
\newcommand{\di}{\,\mathrm{d}}

\renewcommand{\phi}{\varphi}
\renewcommand{\rho}{\varrho}
\renewcommand{\theta}{\vartheta}

\DeclareMathOperator{\Divv}{div}
\renewcommand{\div}{\Divv}

\DeclarePairedDelimiter{\set}{\{}{\}}

\allowdisplaybreaks

\begin{document}

\title[On the monotonicity of weighted perimeters of convex bodies]
{On the monotonicity of weighted perimeters \\ of convex bodies}

\author[G.~Saracco]{Giorgio Saracco}
\address[G.~Saracco]{Dipartimento di Matematica e Informatica ``Ulisse Dini'', Universit\`a di Firenze, viale Morgagni 67/A, 50134 Firenze (FI), Italy}
\email{giorgio.saracco@unifi.it}

\author[G.~Stefani]{Giorgio Stefani}
\address[G.~Stefani]{Scuola Internazionale Superiore di Studi Avanzati (SISSA), via Bonomea 265, 34136 Trieste (TS), Italy}
\email{gstefani@sissa.it {\normalfont or} giorgio.stefani.math@gmail.com}

\date{\today}



\thanks{\textit{Acknowledgements}.
The authors thank Flavia Giannetti for her precious comments on a preliminary version of this work.
The authors are members of the Istituto Nazionale di Alta Matematica (INdAM), Gruppo Nazionale per l'Analisi Matematica, la Probabilit\`a e le loro Applicazioni (GNAMPA).
The first-named author has received funding from INdAM under the INdAM--GNAMPA Project 2023 \textit{Esistenza e propriet\`a fini di forme ottime}, codice CUP\_E53\-C22\-00\-19\-30\-001,
and from Universit\`a di Trento (UNITN) under the Starting Grant Giovani Ricercatori 2021 project \textit{WeiCAp}, codice CUP\_E65\-F21\-00\-41\-60\-001.
The second-named author has received funding from INdAM under the INdAM--GNAMPA 2023 Project \textit{Problemi variazionali per funzionali e operatori non-locali}, codice CUP\_E53\-C22\-00\-19\-30\-001, and from the European Research Council (ERC) under the European Union Horizon 2020 research and innovation program (grant agreement No.~945655).
The present research started during a visit of the first-named author at the Scuola Internazionale Superiore di Studi Avanzati (SISSA). The first-named author wishes to thank SISSA for its kind hospitality.}

\begin{abstract}
We prove that, among weighted isotropic perimeters, only constant multiples of the Euclidean perimeter  satisfy the monotonicity property on nested convex bodies.
Although the analogous result fails for general weighted anisotropic perimeters, a similar characterization holds for radially-weighted   anisotropic densities.
\end{abstract}
 \hspace{-2cm}
 {
 \begin{minipage}[t]{0.6\linewidth}
 \begin{scriptsize}
 \vspace{-3cm}
 This is a pre-print of an article published in \emph{Math.\ Nachr.}. The final authenticated version is available online at: \href{https://doi.org/10.1002/mana.202300280}{https://doi.org/10.1002/mana.202300280}
 \end{scriptsize}
\end{minipage} 
}

\maketitle

\section{Introduction}

\subsection{Monotonicity property}
Let $N\ge2$. If $A,B\subset\R^N$ are two nested convex bodies, that is compact convex sets with non-empty interior such that $A\subset B$, then
\begin{equation}
\label{eq:monotonicity}
P(A)\le P(B),
\end{equation}
where $P(E)=\mathscr H^{N-1}(\partial E)$ denotes the Euclidean perimeter of the convex body $E\subset\R^N$.
The monotonicity property~\eqref{eq:monotonicity} is well known and dates back to the ancient Greeks (Archimedes took it as a postulate in his work on the sphere and the cylinder~\cite{Arc04book}*{p.~36}). 

Inequality~\eqref{eq:monotonicity} can be proved in several ways: by the \emph{Cauchy formula} for the area surface of convex bodies~\cite{BF87book}*{\S7}; by the monotonicity property of \emph{mixed volumes}~\cite{BF87book}*{\S8}; by the Lipschitz property of the projection on a convex closed set~\cite{BFK95}*{Lem.~2.4}; by the fact that the perimeter is decreased under intersection with half-spaces~\cite{Mag12book}*{Ex.~15.13}.  

Inequality~\eqref{eq:monotonicity} extends to the \emph{anisotropic} (\emph{Wulff}) \emph{$\Phi$-perimeter}
\begin{equation*}
P_\Phi(E)=\int_{\partial  E} \Phi\big(\nu_E(x)\big)\di\mathscr H^{N-1}(x),
\end{equation*}   
where $\nu_E\colon\partial E\to\mathbb S^{N-1}$ is the inner unit normal of the convex body $E\subset\R^N$ (defined $\mathscr H^{N-1}$-a.e.\ on $\partial E$) and $\Phi\colon\R^N\to[0,+\infty]$ is a fixed lower-semicontinuous, positively $1$-homogeneous and convex function.
Clearly, if $\Phi=|\cdot|$, then $P_\Phi(E)=P(E)$. 
Similarly to~\eqref{eq:monotonicity}, the monotonicity of the $\Phi$-perimeter is a consequence of one of the following: the \emph{Cauchy formula} for the anisotropic perimeter~\cite{BF87book}*{\S7}; the monotonicity property of \emph{mixed volumes}~\cite{BF87book}*{\S8}; the fact that the anisotropic perimeter is decreased under intersection with half-spaces~\cite{Mag12book}*{Rem.~20.3}.

In passing, we mention that the monotonicity property holds even for perimeter functionals of the \emph{non-local} type, as the \emph{fractional perimeter}~\cite{FFMMM15}*{Lem.~B.1} and, more generally, \emph{non-local perimeters} induced by a suitable interaction kernel~\cite{BS22}*{Cor.~2.30}.

The monotonicity property of perimeters has gained increasing attention in recent years. We refer to~\cites{CGL15,CGLP16,LCL08,Ste18} and to the survey~\cite{G17} for quantitative versions of the monotonicity inequality (see also~\cite{GS2x} for the quantitative monotonicity in the \emph{non-local} setting), and to~\cites{GS23,CGLP19,H23,M22,DLW22,B21,DV20} for some applications and related results.

\subsection{Main result}

In this note, we are interested in studying the monotonicity property on nested convex bodies for the class of \emph{weighted perimeters}.
Given a Borel function $f\colon\R^N\to[0,+\infty]$, we let
\begin{equation}
\label{eq:wpi}
P_f(E)=\int_{\partial E}f(x)\di \mathscr H^{N-1}(x)
\end{equation}
be the \emph{weighted} (\emph{isotropic}) \emph{perimeter} of the convex body $E\subset\R^N$.
Clearly, if $f\equiv c$ for some $c\in[0,+\infty)$, then $P_f=c\,P$, a constant multiple of the Euclidean perimeter. 
Weighted perimeters have been largely investigated in relation to \emph{isoperimetric}, \emph{cluster} and \emph{Cheeger problems}, see~\cites{BCSu,CP17a, CP17b,DPFP17, FPSS22, FPS23a,FPS23b,PS18,PS20, Sar18} and the survey~\cite{P15} for an account on the existing literature.

Our main result is the following rigidity property, namely, the only weighted perimeter satisfying the monotonicity property is (a constant multiple of) the Euclidean perimeter.

\begin{theorem}
\label{thm:main}
Let $f\colon\R^N \to[0,+\infty]$ be a Borel function such that $f\in L^1_{\rm loc}(\R^N)$.
If the weighted perimeter~$P_f$ in~\eqref{eq:wpi} satisfies the monotonicity property, i.e., 
\begin{equation}
\label{eq:monotonicity_w}
\text{
$P_f(A)\le P_f(B)$ for any two nested convex bodies $A\subset B$ in~$\R^N$,
}
\end{equation}
then $f\equiv c$ a.e.\ for some $c\ge0$.
\end{theorem}

\cref{thm:main} is quite intuitive. 
In fact, one clearly expects that, if~$f$ is not constant in some direction, then the monotonicity property should be violated on any suitable family of convex bodies with some side (continuously) deforming along that direction.
However, one should carefully keep into account the values of~$f$ on the entire boundary of each convex body of the family, which forces one to consider deformations in that direction given by graphs of concave functions fixing the boundary of the chosen side. 

One may wonder whether the analog of \cref{thm:main} holds for weighted \emph{anisotropic} perimeters.
More precisely,  given a non-negative Finslerian weight $f\colon\R^N\times\mathbb S^{N-1}\to[0,+\infty]$ (i.e., possibly depending  also on the inner unit normal $\nu_E\colon\partial E\to\mathbb S^{N-1}$ of the convex body $E\subset\R^N$) and assuming the monotonicity of the weighted anisotropic perimeter $P_f$, is it true that $f=f(x,\nu)$ does not depend on~$x$? 
This is in general false.
As a counterexample, consider any bounded vector field $F\in C^1(\R^N;\R^N)$ with constant divergence, $\div F\equiv \alpha$ for some $\alpha\in[0,+\infty)$, and define the anisotropic weight $f\colon\R^N\times\mathbb S^{N-1}\to[0,+\infty)$ as
\begin{equation}
\label{eq:counterex}
f(x,\nu)=F(x)\cdot\nu+\beta
\quad
\text{for}\ x\in\R^N\ \text{and}\ \nu\in\mathbb S^{N-1},
\end{equation}
where $\beta \in [\|F\|_\infty,+\infty)$ ensures the non-negativity of the weight $f$. By the Divergence Theorem, the anisotropic weighted perimeter
\begin{equation}
\label{eq:wpa}
P_f(E)=\int_{\partial E}f(x,\nu_E(x))\di\mathscr H^{N-1}(x)
\end{equation}
on the convex body $E\subset\R^N$ satisfies
\begin{align*}
P_f(E) 
&
= \int_{\partial E} f(x,\nu_E(x)) \di\mathscr H^{N-1}(x) 
\\
&
= \int_{\partial E} F(x)\cdot\nu_E(x)\di\mathscr H^{N-1}(x) + \beta\, P(E)
\\
&
= 
\int_E \div F \di x + \beta\, P(E)
\\
& = \alpha\,|E|+\beta\, P(E),
\end{align*}
readily yielding the desired monotonicity property in virtue of that of the Euclidean perimeter~\eqref{eq:monotonicity} and that of the Lebesgue measure with respect to nestedness.

Despite the counterexample in~\eqref{eq:counterex}, from \cref{thm:main} we can deduce the following result, which provides a partial analog of the rigidity property in the anisotropic regime under some additional structural assumptions on the weight function.

\begin{corollary}
\label{res:radial}
Let $f\colon\R^N\times\R^N\to[0,+\infty]$ be a Borel function such that $f\in L^1_{\rm loc}(\R^{2N})$.
Assume that there exist a radial Borel function $g\colon\R^N\to[0,+\infty]$ and a lower semicontinuous, $1$-homogeneous and convex function $\Phi\colon\R^N\to(0,+\infty]$ such that
\begin{equation}
\label{eq:radial_ass}
f(x,v)=g(x)\,\Phi(v)
\quad 
\text{for}\ x\in\R^N\ \text{and}\ v\in\R^N.
\end{equation}
If the anisotropic weighted perimeter~$P_f$ in~\eqref{eq:wpa} satisfies the monotonicity property~\eqref{eq:monotonicity_w}, then $g\equiv c$ a.e.\ for some $c\ge0$.
\end{corollary}

The proof of \cref{res:radial} combines the invariance of the monotonicity property with respect to rotations with \cref{thm:main}.
%We believe that the radial symmetry of the function $g$ in the hypotheses of \cref{res:radial} can be dropped, but we leave this possible improvement to future works.

\section{Proofs of the statements}
\label{sec:proof}

\subsection{Proof of \texorpdfstring{\cref{thm:main}}{Theorem 1.1}}

We begin by observing that it is not restrictive to assume that $f\in C^\infty(\R^N)$.
Indeed, given $A\subset B$ two nested convex bodies in~$\R^N$,  the translated sets $A+y\subset B+y$ are still two nested convex bodies for  any $y\in\R^N$. Therefore, in virtue of~\eqref{eq:monotonicity_w} and changing variables, we get
\begin{equation}
\label{eq:translation}
\int_{\partial A}f(x-y)\di\mathscr H^{N-1}(x)
\le 
\int_{\partial B}f(x-y)\di\mathscr H^{N-1}(x).
\end{equation} 
Let now $(\rho_\eps)_{\eps>0}\subset C^\infty_c(\R^N)$ be any family of non-negative convolution kernels (for instance, $\rho_\eps=\eps^{-N}\rho(\sfrac\cdot\eps)$ for some $\rho\in C^\infty(\R^N)$ such that $\operatorname{supp}\rho\subset B_1$, $\rho\ge0$, and $\int_{\R^N}\rho\di x=1$). 
Multiplying~\eqref{eq:translation} by $\rho_\eps(y)$, integrating on $\R^N$ with respect to $y$, and owing to the Fubini--Tonelli Theorem, we infer that
\begin{align*}
\int_{\partial A} f_\eps(x)\di\mathscr H^{N-1}(x) 
&=
\int_{\partial A}\int_{\R^N} f(x-y)\,\rho_\eps(y)\di y\di\mathscr H^{N-1}(x)
\\
&\le 
\int_{\partial B}\int_{\R^N} f(x-y)\,\rho_\eps(y)\di y\di\mathscr H^{N-1}(x)
=
\int_{\partial B} f_\eps(x)\di\mathscr H^{N-1}(x)
,
\end{align*}
where $f_\eps=f*\rho_\eps\in C^\infty(\R^N)$ is the standard convolution. By the arbitrariness of the nested convex bodies $A$ and $B$, the weight $f_\eps$ still verifies~\eqref{eq:monotonicity_w} for each $\eps>0$. If we show that $\nabla f_\eps\equiv0$ for each $\eps>0$, then also $\nabla f\equiv0$ in the sense of distributions, and thus $f$ is equivalent to a constant function.

Consequently, from now on, we assume that $f\in C^\infty(\R^N)$.
We now claim that $\partial_{x_N}f(x)=0$ for each $x\in\R^N$. 
By the translation invariance in~\eqref{eq:translation}, we just need to show that $\partial_{x_N}f(0)=0$.

\begin{figure}[t]
\begin{tikzpicture}
\begin{axis}[width=200mm,
            height=200mm,
            axis lines=none,
            line width=1pt,
            xmin=-5.5,
            xmax=5.5,
            ymin=-5.5,
            ymax=5.5,
            zmax=2,
            zmin=-4]
 
\addplot3 [
    domain=-2:2,
    domain y = -2:2,
    samples = 40,
    samples y = 40,
    surf,
    fill=gray,
    faceted color=gray,
    opacity=0.1,
    fill opacity=0.3] {0.3*((x^2-4)*(y^2-4))^(1/2)};
	\draw[dashed] (-2,-2,0) -- (-2,2,0) -- (2,2,0);
        \draw (2,2,0) -- (2,-2,0) -- (-2,-2,0);
	\draw[dashed] (-2,-2,-1) -- (-2,2,-1) -- (2,2,-1);
        \draw (2,2,-1) -- (2,-2,-1) -- (-2,-2,-1);
        \fill[fill=black] (0,0,0) circle (1pt); 
        \draw[-latex] (0,0,0) -- node[left]{$e_3$} (0,0,0.5) ;
	\draw (-2,-2,0) -- (-2,-2,-1);
	\draw[dashed] (-2,2,0) -- (-2,2,-1);
	\draw (2,2,0) -- (2,2,-1);
	\draw (2,-2,0) -- (2,-2,-1);
    \node at (0,0,-1.25) {$E(\lambda)$};
    \node at (1,1.5,1) {$\Gamma_h(\lambda)$};
\end{axis}
\end{tikzpicture}
\caption{The set $E(\lambda)$ and its deformation $F_h(\lambda)$ for $\lambda>0$ and a given concave function $h\colon[-\delta, \delta]^{2}\to\R$ vanishing on the boundary of its domain.\label{fig:sets}}
\end{figure}
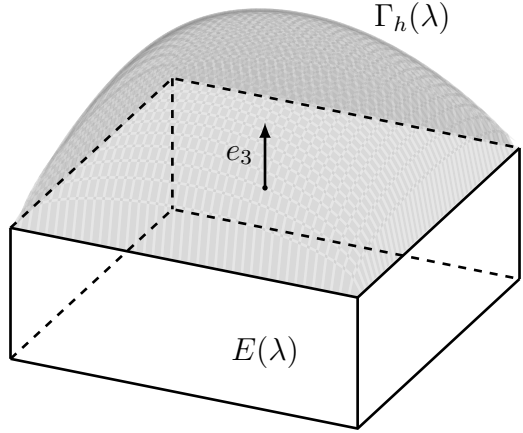

Let $\delta>0$ to be chosen later on. For $\lambda\in\R$, we define
\begin{equation*}
E(\lambda)=
\begin{cases}
[-\delta,\delta]^{N-1}\times [-\delta, 0] 
& 
\text{for}\ \lambda\ge0,
\\[1ex]
[-\delta,\delta]^{N-1}\times [0,\delta] 
&
\text{for}\ \lambda<0.
\end{cases}
\end{equation*}
Moreover, given $h\colon[-\delta, \delta]^{N-1}\to\R$ any concave function vanishing on the boundary of $[-\delta, \delta]^{N-1}\subset\R^{N-1}$, we set
\begin{equation*}
\Gamma_h(\lambda)=
\begin{cases}
\set*{x=(x',x_N)\in\R^N : x'\in[-\delta,\delta]^{N-1}\ \text{and}\ 0\le x_N\le \lambda h(x')}
& 
\text{for}\ \lambda\ge0,
\\[1.5ex]
\set*{x=(x',x_N)\in\R^N : x'\in[-\delta,\delta]^{N-1}\ \text{and}\ \lambda h(x')\le x_N\le 0}&
\text{for}\ \lambda<0,
\end{cases}
\end{equation*}
and we refer to \cref{fig:sets} for a visual aid in the $3$-dimensional case.
%%%%
%%%%%
Note that $E(\lambda)$ and $F(\lambda)=E(\lambda)\cup \Gamma_h(\lambda)$ are convex bodies in $\R^N$ with $E(\lambda)\subset F(\lambda)$ for all $\lambda\in\R$. Hence, in virtue of~\eqref{eq:monotonicity_w}, we get that
\begin{equation*}
P_f(F(\lambda)) \ge P_f(E(\lambda))
\quad
\text{for all}\ \lambda\in\R.
\end{equation*}
By the area formula, the above inequality rewrites as 
\begin{equation*}
\int_{[-\delta,\delta]^{N-1}} f(x', \lambda h(x'))\, \sqrt{1+\lambda^2\,|\nabla h(x')|^2} - f(x', 0)\di x'\ge 0
\quad
\text{for all}\ \lambda\in\R.
\end{equation*} 
In particular, the function $\ell\colon\R\to[0,+\infty)$, given by 
\begin{equation*}
\ell(\lambda)
=
\int_{[-\delta,\delta]^{N-1}} f(x', \lambda h(x'))\, \sqrt{1+\lambda^2\,|\nabla h(x')|^2}\di x'
\quad
\text{for}\ \lambda\in\R,
\end{equation*}
achieves its minimum at $\lambda=0$, so that $\ell'(0)=0$.
Since $f$ is a smooth function, we can exchange the differentiation and the integration signs, obtaining
\begin{equation}
\label{eq:EL}
\int_{[-\delta,\delta]^{N-1}} \partial_{x_N}f(x',0)\,h(x')\di x'
=0
\end{equation}
for any concave function $h\colon[-\delta, \delta]^{N-1}\to\R$ vanishing on the boundary of $[-\delta, \delta]^{N-1}$.
In particular, $h(x')>0$ for all $x'\in(-\delta,\delta)^{N-1}$ as soon as $h\not\equiv0$.
By contradiction, if $\partial_{x_N}f(0)\ne0$, then, by smoothness of~$f$, we may assume that $\partial_{x_N} f(x)$ has constant sign for each $x\in B_r(0)$ for some $r>0$.
Choosing $\delta>0$ so small that $[-\delta,\delta]^{N-1}\times\set*{0}\subset B_r(0)$, the equality~\eqref{eq:EL} immediately yields a contradiction.

In the previous argument, the choice of fixing the $N$-th component does not play any role and can be repeated almost \textit{verbatim} to show that $\partial_{x_i}f(0)=0$ for each $i=1,\dots,N$. Thus, again by the translation invariance~\eqref{eq:translation}, we get that $\nabla f(x)=0$ for all $x\in\R^N$, yielding the conclusion.
\qed

\begin{remark}
In the above proof, one needs much less than the monotonicity of the perimeter on nested convex bodies in order to conclude that the weight is constant. Indeed, it would be enough to know that, for each direction $\mathrm{e}_i\in\mathbb S^{N-1}$, $i=1,\dots,N$, and each point $x\in\R^N$, the monotonicity property holds on two  hypercubes (not necessarily with the same edge size) with a face containing $x$ and orthogonal to $\mathrm{e}_i$ with opposite outward normals $\pm\mathrm{e}_i$ on that face.    
\end{remark}

\subsection{Proof of \texorpdfstring{\cref{res:radial}}{Corollary 1.2}}

Let us denote by $\mathrm{SO}(N)$ be the special orthogonal group, and let $\mu\in\mathscr P(\mathrm{SO}(N))$ be the (unique) \emph{Haar} \emph{probability} \emph{measure} on $\mathrm{SO}(N)$ (see~\cite{DS14book} for a detailed exposition).
Given $A\subset B$ two nested convex bodies in~$\R^N$, the rotated sets $\mathcal R (A)\subset\mathcal R(B)$ are still two nested convex bodies for any $\mathcal R\in\mathrm{SO}(N)$.
Therefore, in virtue of~\eqref{eq:monotonicity_w} and changing variables, we get
\begin{equation}
\label{eq:rotation}
\int_{\partial A} f\big(\mathcal R(x),\mathcal R(\nu_A(x))\big)\di\mathscr H^{N-1}(x)
\le
\int_{\partial B} f\big(\mathcal R(x),\mathcal R(\nu_B(x))\big)\di\mathscr H^{N-1}(x),
\end{equation}
owing to the elementary facts that $\mathcal R(\partial E) = \partial \mathcal R(E)$ and that $\nu_{\mathcal R(E)}(\mathcal R(x))=\mathcal R(\nu_E(x))$ for $\mathscr H^{N-1}$-a.e.\  $x\in\partial E$ whenever $E\subset\R^N$ is a convex body (refer to~\cite{Mag12book}*{\S17.1} for a precise justification).
Due to~\eqref{eq:radial_ass} and the radial assumption on $g$, inequality~\eqref{eq:rotation} rewrites as
\begin{equation}
\label{eq:rotation_radial}
\int_{\partial A} g(x)\,\Phi\big(\mathcal R(\nu_A(x))\big)\di\mathscr H^{N-1}(x)
\le
\int_{\partial B} g(x)\,\Phi\big(\mathcal R(\nu_B(x))\big)\di\mathscr H^{N-1}(x)
\end{equation}
for $\mathcal R\in\mathrm{SO}(N)$.
We now claim that the function
\begin{equation*}
\mathbb S^{N-1}\ni\nu\mapsto
\int_{\mathrm{SO}(N)}\Phi\big(\mathcal R(\nu)\big)\di\mu(\mathcal R)
\end{equation*}
is constant. 
Indeed, given any $\nu\in\mathbb S^{N-1}$, we can find $\mathcal R_\nu\in\mathrm{SO}(N)$ such that $\nu=\mathcal R_\nu(\mathrm{e}_1)$. 
Due to the invariance properties of the Haar measure~$\mu$, we can compute
\begin{align}
\int_{\mathrm{SO}(N)}\Phi\big(\mathcal R(\nu)\big)\di\mu(\mathcal R)
\nonumber
&=
\int_{\mathrm{SO}(N)}\Phi\big(\mathcal R(\mathcal R_\nu(\mathrm{e}_1))\big)\di\mu(\mathcal R)
\nonumber
\\
&=
\int_{\mathrm{SO}(N)}\Phi\big(\mathcal Q(\mathrm{e}_1)\big)\di\mu(\mathcal Q\mathcal R^{-1}_\nu)
\nonumber
\\
&=
\int_{\mathrm{SO}(N)}\Phi\big(\mathcal Q(\mathrm{e}_1)\big)\di\mu(\mathcal Q)
\label{eq:int_nonzero}
\end{align} 
where, with a slight abuse of notation, $\mathcal Q\mapsto\mu(\mathcal Q\mathcal R^{-1}_\nu)$ stands for the push-forward of the measure $\mu$ with respect to the right translation by $\mathcal R_\nu^{-1}$. 
Hence, integrating on $\mathrm{SO}(N)$ with respect to $\mu$, using the Fubini--Tonelli Theorem, the above equality, that $\Phi>0$, and simplifying, from~\eqref{eq:rotation_radial} we get
\begin{equation*}
\int_{\partial A} g(x)
\di\mathscr H^{N-1}(x)
\le
\int_{\partial B} g(x)
\di\mathscr H^{N-1}(x)
\end{equation*}
for any two nested convex bodies $A\subset B$. 
The conclusion follows from \cref{thm:main}.
\qed

\begin{remark}
One could slightly weaken the hypotheses of \cref{res:radial} by allowing $\Phi$ to also attain zero. 
In fact, it is enough to require that the integral in~\eqref{eq:int_nonzero} is not zero.
\end{remark}

%%%  BIBLIO %%%


\begin{bibdiv}
\begin{biblist}

\bib{Arc04book}{book}{
      author={Archimedes},
       title={The {W}orks of {A}rchimedes. {V}ol. {I}},
   publisher={Cambridge University Press, Cambridge},
        date={2004},
        ISBN={0-521-66160-9},
%         url={https://doi.org/10.1017/CBO9780511482557},
%        note={The two books on the sphere and the cylinder, Translated into English, together with Eutocius' commentaries, with commentary, and critical edition of the diagrams by Reviel Netz},
      review={\MR{2093668}},
}

\bib{BCSu}{article}{
      author={Beck, Lisa},
      author={Cinti, Eleonora},
      author={Seis, Christian},
       title={Optimal regularity of isoperimetric sets with {H}\"older
  densities},
  journal={Calc. Var. Partial Differential Equations},
   volume={62},
   date={2023},
   pages={214},
}

\bib{BS22}{article}{
   author={Bessas, Konstantinos},
   author={Stefani, Giorgio},
   title={Non-local BV functions and a denoising model with $L^1$ fidelity},
   journal={Adv. Calc. Var.},
   date={to appear},
   note={Preprint available at \href{https://arxiv.org/abs/2210.11958}{arXiv:2210.11958}},
}

\bib{B21}{article}{
   author={Berman, Robert J.},
   title={Convergence rates for discretized Monge--Amp\`ere equations and
   quantitative stability of optimal transport},
   journal={Found. Comput. Math.},
   volume={21},
   date={2021},
   number={4},
   pages={1099--1140},
   issn={1615-3375},
   review={\MR{4298241}},
%   doi={10.1007/s10208-020-09480-x},
}

\bib{BF87book}{book}{
      author={Bonnesen, T.},
      author={Fenchel, W.},
       title={Theory of {C}onvex {B}odies},
   publisher={BCS Associates, Moscow, ID},
        date={1987},
        ISBN={0-914351-02-8},
        note={Translated from the German and edited by L. Boron, C. Christenson
  and B. Smith},
      review={\MR{920366}},
}

\bib{BFK95}{article}{
      author={Buttazzo, Giuseppe},
      author={Ferone, Vincenzo},
      author={Kawohl, Bernhard},
       title={Minimum problems over sets of concave functions and related
  questions},
        date={1995},
        ISSN={0025-584X},
     journal={Math. Nachr.},
      volume={173},
       pages={71\ndash 89},
         url={https://doi.org/10.1002/mana.19951730106},
      review={\MR{1336954}},
}

\bib{CGL15}{article}{
      author={Carozza, Menita},
      author={Giannetti, Flavia},
      author={Leonetti, Francesco},
      author={Passarelli di Napoli, Antonia},
       title={A sharp quantitative estimate for the perimeters of convex sets
  in the plane},
        date={2015},
        ISSN={0944-6532},
     journal={J. Convex Anal.},
      volume={22},
      number={3},
       pages={853\ndash 858},
      review={\MR{3400158}},
}

\bib{CGLP16}{article}{
      author={Carozza, Menita},
      author={Giannetti, Flavia},
      author={Leonetti, Francesco},
      author={Passarelli di Napoli, Antonia},
       title={A sharp quantitative estimate for the surface areas of convex
  sets in {$\mathbb{R}^3$}},
        date={2016},
        ISSN={1120-6330},
     journal={Atti Accad. Naz. Lincei Rend. Lincei Mat. Appl.},
      volume={27},
      number={3},
       pages={327\ndash 333},
         url={https://doi.org/10.4171/RLM/737},
      review={\MR{3510903}},
}

\bib{CGLP19}{article}{
   author={Carozza, Menita},
   author={Giannetti, Flavia},
   author={Leonetti, Francesco},
   author={Passarelli di Napoli, Antonia},
   title={Convex components},
   journal={Commun. Contemp. Math.},
   volume={21},
   date={2019},
   number={6},
   pages={1850036, 10},
   issn={0219-1997},
   review={\MR{3996972}},
%   doi={10.1142/S0219199718500360},
}

\bib{CP17a}{article}{
      author={Cinti, E.},
      author={Pratelli, A.},
       title={Regularity of isoperimetric sets in {$\mathbb{R}^2$} with
  density},
        date={2017},
        ISSN={0025-5831},
     journal={Math. Ann.},
      volume={368},
      number={1-2},
       pages={419\ndash 432},
         url={https://doi.org/10.1007/s00208-016-1409-y},
      review={\MR{3651579}},
}

\bib{CP17b}{article}{
      author={Cinti, Eleonora},
      author={Pratelli, Aldo},
       title={The {$\varepsilon-\varepsilon^\beta$} property, the boundedness
  of isoperimetric sets in {$\mathbb{R}^N$} with density, and some
  applications},
        date={2017},
        ISSN={0075-4102},
     journal={J. Reine Angew. Math.},
      volume={728},
       pages={65\ndash 103},
         url={https://doi.org/10.1515/crelle-2014-0120},
      review={\MR{3668991}},
}

\bib{DV20}{article}{
   author={Dekeyser, Justin},
   author={Van Schaftingen, Jean},
   title={Range convergence monotonicity for vector measures and range
   monotonicity of the mass},
   journal={Ric. Mat.},
   volume={69},
   date={2020},
   number={1},
   pages={293--326},
   issn={0035-5038},
   review={\MR{4098186}},
%   doi={10.1007/s11587-019-00463-x},
}

\bib{DPFP17}{article}{
      author={De~Philippis, Guido},
      author={Franzina, Giovanni},
      author={Pratelli, Aldo},
       title={Existence of isoperimetric sets with densities ``converging from
  below'' on {$\mathbb{R}^N$}},
        date={2017},
        ISSN={1050-6926},
     journal={J. Geom. Anal.},
      volume={27},
      number={2},
       pages={1086\ndash 1105},
         url={https://doi.org/10.1007/s12220-016-9711-1},
      review={\MR{3625144}},
}

\bib{DS14book}{book}{
   author={Diestel, Joe},
   author={Spalsbury, Angela},
   title={The {J}oys of {H}aar {M}easure},
   series={Graduate Studies in Mathematics},
   volume={150},
   publisher={American Mathematical Society, Providence, RI},
   date={2014},
   pages={xiv+320},
   isbn={978-1-4704-0935-7},
   review={\MR{3186070}},
%   doi={10.1090/gsm/150},
}

\bib{DLW22}{article}{
   author={Du, Qiang},
   author={Lu, Xin Yang},
   author={Wang, Chong},
   title={The average-distance problem with an Euler elastica penalization},
   journal={Interfaces Free Bound.},
   volume={24},
   date={2022},
   number={1},
   pages={137--162},
   issn={1463-9963},
   review={\MR{4395706}},
%   doi={10.4171/ifb/470},
}

\bib{FFMMM15}{article}{
   author={Figalli, A.},
   author={Fusco, N.},
   author={Maggi, F.},
   author={Millot, V.},
   author={Morini, M.},
   title={Isoperimetry and stability properties of balls with respect to
   nonlocal energies},
   journal={Commun. Math. Phys.},
   volume={336},
   date={2015},
   number={1},
   pages={441--507},
   issn={0010-3616},
   review={\MR{3322379}},
%   doi={10.1007/s00220-014-2244-1},
}

\bib{FPSS22}{article}{
  author = {Franceschi, Valentina},
  author = {Pinamonti, Andrea},
  author = {Saracco, Giorgio},
  author = {Stefani, Giorgio },
  title = {The Cheeger problem in abstract measure spaces},
  journal = {J. London Math. Soc.},
  year = {to appear},
  note={Preprint available at \href{https://arxiv.org/abs/2207.00482}{arXiv:2207.00482}},
}

\bib{FPS23b}{article}{
      author={Franceschi, Valentina},
      author={Pratelli, Aldo},
      author={Stefani, Giorgio},
       title={On the {S}teiner property for planar minimizing clusters. {T}he
  anisotropic case.},
        date={2023},
     journal={J. {\'E}c. Polytech. Math.},
      volume={10},
       pages={989\ndash 1045},
       review={\MR{4600399}},
}

\bib{FPS23a}{article}{
      author={Franceschi, Valentina},
      author={Pratelli, Aldo},
      author={Stefani, Giorgio},
       title={On the {S}teiner property for planar minimizing clusters. {T}he
  isotropic case},
        date={2023},
        ISSN={0219-1997},
     journal={Commun. Contemp. Math.},
      volume={25},
      number={5},
       pages={Paper No. 2250040, 29},
%         url={https://doi.org/10.1142/S0219199722500407},
      review={\MR{4579987}},
}

\bib{G17}{article}{
   author={Giannetti, Flavia},
   title={Sharp geometric quantitative estimates},
   journal={Atti Accad. Naz. Lincei Rend. Lincei Mat. Appl.},
   volume={28},
   date={2017},
   number={1},
   pages={1--6},
   issn={1120-6330},
   review={\MR{3621767}},
%   doi={10.4171/RLM/748},
}

\bib{GS23}{article}{
   author={Giannetti, Flavia},
   author={Stefani, Giorgio},
   title={On the convex components of a set in $\mathbb{R}^{n}$},
   journal={Forum Math.},
   volume={35},
   date={2023},
   number={1},
   pages={187--199},
   issn={0933-7741},
   review={\MR{4529426}},
%   doi={10.1515/forum-2022-0203},
}

\bib{GS2x}{article}{
   author={Giannetti, Flavia},
   author={Stefani, Giorgio},
   title={On the monotonicity of non-local perimeter of convex bodies},
   year={2023}
   note={Preprint available at \href{https://arxiv.org/abs/2309.11296}{arXiv:2309.11296}},
}

\bib{H23}{article}{
   author={Hynd, Ryan},
   title={A doubly monotone flow for constant width bodies in $\mathbb{R}^3$},
   conference={
      title={Geometric and Functional Inequalities and Recent Topics in
      Nonlinear PDEs},
   },
   book={
      series={Contemp. Math.},
      volume={781},
      publisher={Amer. Math. Soc., Providence, RI},
   },
   date={2023},
   pages={49--101},
   review={\MR{4531939}},
%   doi={10.1090/conm/781/15709},
}

\bib{LCL08}{article}{
      author={La~Civita, Marianna},
      author={Leonetti, Francesco},
       title={Convex components of a set and the measure of its boundary},
        date={2008/09},
        ISSN={1825-1269},
     journal={Atti Semin. Mat. Fis. Univ. Modena Reggio Emilia},
      volume={56},
       pages={71\ndash 78},
      review={\MR{2604730}},
}

\bib{Mag12book}{book}{
      author={Maggi, Francesco},
       title={Sets of {F}inite {P}erimeter and {G}eometric {V}ariational {P}roblems},
      series={Cambridge Studies in Advanced Mathematics},
   publisher={Cambridge University Press, Cambridge},
        date={2012},
      volume={135},
        ISBN={978-1-107-02103-7},
%         url={https://doi.org/10.1017/CBO9781139108133},
%        note={An introduction to geometric measure theory},
      review={\MR{2976521}},
}

\bib{M22}{article}{
   author={Melchionna, Andrew},
   title={The sandpile identity element on an ellipse},
   journal={Discrete Contin. Dyn. Syst.},
   volume={42},
   date={2022},
   number={8},
   pages={3709--3732},
   issn={1078-0947},
   review={\MR{4447555}},
%   doi={10.3934/dcds.2022029},
}

\bib{P15}{article}{
   author={Pratelli, Aldo},
   title={A survey on the existence of isoperimetric sets in the space
   $\mathbb{R}^N$ with density},
   journal={Atti Accad. Naz. Lincei Rend. Lincei Mat. Appl.},
   volume={26},
   date={2015},
   number={1},
   pages={99--118},
   issn={1120-6330},
   review={\MR{3345326}},
%   doi={10.4171/RLM/696},
}

\bib{PS18}{article}{
      author={Pratelli, Aldo},
      author={Saracco, Giorgio},
       title={On the isoperimetric problem with double density},
        date={2018},
        ISSN={0362-546X},
     journal={Nonlinear Anal.},
      volume={177},
      number={part B},
       pages={733\ndash 752},
%         url={https://doi.org/10.1016/j.na.2018.04.009},
      review={\MR{3886599}},
}

\bib{PS20}{article}{
      author={Pratelli, Aldo},
      author={Saracco, Giorgio},
       title={The {$\varepsilon-\varepsilon^{\beta}$} property in the
  isoperimetric problem with double density, and the regularity of
  isoperimetric sets},
        date={2020},
        ISSN={1536-1365},
     journal={Adv. Nonlinear Stud.},
      volume={20},
      number={3},
       pages={539\ndash 555},
%         url={https://doi.org/10.1515/ans-2020-2074},
      review={\MR{4129341}},
}

\bib{Sar18}{article}{
    AUTHOR = {Saracco, Giorgio},
     TITLE = {Weighted {C}heeger sets are domains of isoperimetry},
   JOURNAL = {Manuscripta Math.},
    VOLUME = {156},
      YEAR = {2018},
    NUMBER = {3-4},
     PAGES = {371\ndash381},
      ISSN = {0025-2611},
  review = {\MR{3811794}},
       %DOI = {10.1007/s00229-017-0974-z},
      % URL = {https://doi.org/10.1007/s00229-017-0974-z},
}

\bib{Ste18}{article}{
      author={Stefani, Giorgio},
       title={On the monotonicity of perimeter of convex bodies},
        date={2018},
        ISSN={0944-6532},
     journal={J. Convex Anal.},
      volume={25},
      number={1},
       pages={93\ndash 102},
      review={\MR{3756927}},
}

\end{biblist}
\end{bibdiv}
\end{document}